\magnification=\magstephalf
\input amssym.def
\input amssym.tex

\def\N{{\Bbb N}}
\def\Z{{\Bbb Z}}
\def\Q{{\Bbb Q}}
\def\R{{\Bbb R}}

\def\F{{\Bbb F}}

\def\emptyset{\varnothing}

\def\geq{\geqslant}
\def\leq{\leqslant}

\def\hat{\widefat}
\def\gr{{\rm gr}}

\def\AM/{d'Abfyan\-kar--Mof}

\def\hat{\widehat}
\def\gr{{\rm gr}}


\centerline {\bf CHARACTERISTIC POLYHEDRA OF SINGULARITIES }
\centerline {\bf WITHOUT COMPLETION}
\bigskip \bigskip

\centerline { Vincent COSSART and Olivier PILTANT}

\centerline {\tt cossart@math.uvsq.fr, piltant@math.uvsq.fr}

\bigskip\bigskip
\noindent Laboratoire de Math{\'e}matiques LMV UMR 8100 \hfill

\noindent Universit{\'e} de Versailles \hfill

\noindent 45 avenue des {\'E}tats-Unis\hfill

\noindent 78035 VERSAILLES Cedex (France)
\bigskip

\hfill{}{\it Dedicated to Heisuke Hironaka on his eightieth birthday.}

\bigskip \bigskip

\centerline {\underbar{INTRODUCTION}}
\medskip

\medskip
\noindent {\bf }

Let $(R,M, k:=R/M)$ be a regular local ring, $f\in M$ and $(y,u_1,\dots,u_d)$ be a regular system of parameters (r.s.p. for short) of $R$.
Assume furthermore that
$$
f \not \in (u_1, \ldots ,u_d).
$$

To this situation, H. Hironaka [H] attaches a polyhedron $\Delta(f;u_1,\dots,u_d;y)\subset  \R_{\geq 0}^d$ and shows
([H] theorem 4.8) that there exists some $z\in \hat{R}$ such that $(z,u_1,\dots,u_d)$ is a r.s.p. of $\hat{R}$ and
$$
\Delta(f;u_1,\dots,u_d;z)= \bigcap_{(\hat{y},u_1,\dots,u_d)}\Delta(f;u_1,\dots,u_d;\hat{y}). \eqno (0.1)
$$
Here, the intersection runs over all r.s.p's of $\hat{R}$ of the form $(\hat{y},u_1,\dots,u_d)$.
We prove in theorem {\bf II.3} below that one can choose $z\in R$ whenever $R$ is a G-ring. This is proved
by using a different algorithm than Hironaka's for computing the {\it characteristic} polyhedron
(right hand side of (0.1)).
We emphasize that no assumption is made on the characteristic nor residue characteristic of $R$.
See [C] for the special case of a ring of analytic functions, in a more general setup.

\medskip

Hironaka's polyhedra are an essential tool for defining resolution invariants of singularities. In the introduction
of [H], H. Hironaka explains how to use the characteristic polyhedron for resolving surface singularities ($d=2$).
It is stated how blowing up closed points sharpens the polyhedron and eventually reduces the multiplicity
or produces some permissible curve passing through the singularity.
Thus Resolution of surface Singularities is essentially reduced to a combinatorial statement
about the transformation law for the characteristic polyhedron under blowing up a closed point ([CGO] lemmas 6.2
and 6.3 in the first author's contribution).

\smallskip

In dimension $d\geq 3$, the situation is much more delicate (see comments at the end of the introduction of [H]) and
little is known on birational Resolution of Singularities in positive residue characteristic.
Dimension three Resolution is achieved in [CP1,2] and [CP3], respectively in positive and in mixed characteristic.

\smallskip

Since Hironaka's celebrated equicharacteristic zero Resolution of Singularities, it is expected
that {\it permissible blowing up centers} (i.e. regular centers along which a singular scheme is normally flat)
play a central role in building up a resolution.
Suppose that $A$ is an excellent regular domain, $f \in A$ and let 
$$
{\cal X}:={\rm Spec}(A/(f)).
$$
A typical situation occurs when the generic point $\xi$ of a blowing up center $Z \subset {\cal X}$ has to be blown up.
At a special point $x \in Z$, a characteristic polyhedron $\Delta_x (f;{\bf u}_x)$ is attached
by the formula on the right hand side of (0.1) w.r.t. $R:=A_{\frak{m}}$, $\frak{m}$ the ideal of $x$. 
One is then led to study how the characteristic polyhedron and derived coefficient ideals or numerical invariants 
transform above special points $x\in Z$ when blowing up $Z$. 
This analysis is in principle hopeless if (0.1) requires using nonalgebraic formal coordinates
at each $x \in Z$ for computing $\Delta_x (f;{\bf u}_x)$, as it is expected from [H] theorem 4.8.
Theorem {\bf II.3} then appears as a very useful tool in order to prove some form of coherence or
semicontinuity for resolution invariants derived from the polyhedron function
$$
\Delta: \ x \mapsto \Delta_x (f;{\bf u}_x). \eqno (0.2)
$$
The computation of the $\Delta$-function and application to Resolution of Singularities
is explained in example {\bf II.6} for a concrete arithmetical threefold. An application of theorem {\bf II.3} for
${\cal X}$ a $p$-cyclic covering of any excellent regular germ ${\rm Spec}S$ of residue characteristic $p>0$
is given in [CP3] corollary 3.11: a numerical function $x \mapsto (m(x),\omega (x),\kappa (x))$ refining the multiplicity function
$x \mapsto m(x)$ is defined on ${\cal X}$ and is proved to be a {\it constructible} function.

\medskip

The present result has been exposed by the first author at the conference
``Resolution of singularities and related topics, 80th birthday of Heisuke Hironaka'', Tordesillas.
The authors acknowledge fruitful discussions with U.~Jannsen and S. Saito during the first author's stay at the Regensburg University.

\bigskip 

\centerline {\underbar{I CHARACTERISTIC POLYHEDRON OF $f$}}

\medskip

In this section, we briefly survey known material on $F$-subsets, characteristic polyhedra
and graded algebras with respect to weighted monomial filtrations.

\bigskip

\noindent {\bf I.1 DEFINITION.}{\it

\noindent (i) An $F$-subset $\Delta \subset \R_{\geq 0}^e$ is a closed convex
subset of $\R_{\geq 0}^e$ such that ${\bf v} \in \Delta$ implies ${\bf v}+\R_{\geq 0}^e \subseteq \Delta$;

\noindent (ii) a point ${\bf v}\in \Delta$ is called a vertex if there exists a
positive linear form $L$ on $\R^e$ (i.e. with positive coefficients)
such that
$$\{{\bf v}\}=\Delta\cap\{{\bf a}\in \R_{\geq 0}^e|L({\bf a})=1\};$$

\noindent (iii) an $F$-subset $\Delta$ is called a rational polyhedron if there exists finitely many
nonnegative rational linear forms $L_1, \ldots ,L_n$ on $\R^e$ (i.e. with nonnegative rational coefficients) such that
$$\Delta=\{ {\bf a}\in \R_{\geq 0}^e|L_i({\bf a})\geq 1, \ 1 \leq i \leq n\}.$$
}

\medskip

Given a r.s.p. $(y,u_1,u_2,\dots,u_d)=:(y,{\bf u})$ of a regular local ring $R$ and $f\in R$,
there exists a {\it finite} sum expansion
$$
f=\sum_{{\bf a},b} C_{{\bf a},b}y^b {\bf u}^{\bf a}, \ b\in \N,\ {\bf a}\in \N^d.
\eqno (1.1)$$
where each $C_{{\bf a},b}$ is a unit in $R$. This follows easily from the facts that $R$ is Noetherian
and the map $R \subseteq \hat{R}$ is faithfully flat.

\medskip

Assume furthermore that
$$
f \in M, \ f \not \in (u_1, \ldots ,u_d). \eqno (1.2)
$$

\medskip

We let $\overline{R}:=R/(u_1, \ldots ,u_d)$, $\overline{f}\in \overline{R}$ be the image of $f$ and ``${\rm ord}$''
be the valuation of the discrete valuation ring $\overline{R}$. We let
$$
m:= {\rm ord}\overline{f}\geq 1. \eqno (1.3)
$$

Assumption (1.2) and notation (1.3) are maintained all along this article. We regard ${\bf u}$
as ``fixed" parameters and $y$ as ``varying", which is reflected in the indexing below.

\medskip

For a given expansion (1.1), we denote by ${\rm Supp} (f)\subset \N^{d+1}$ the {\it support} of $f$
in the given expansion, i.e. the set of all $({\bf a},b)$ appearing in (1.1) with some unit coefficient $C_{{\bf a},b}$. Let
$$
{\bf p}: \R^{d+1} \backslash \{b=m\}\longrightarrow \R^d, \ ({\bf a},b)\mapsto {{\bf a} \over m-b}
$$
be the projection on the hyperplane $\{b=0\}$ from the point $({\bf 0},m)$. Define
$$
S(f):={\bf p} \left ({\rm Supp} (f) \cap \{ ({\bf a},b) \in \N^{d+1} , 0\leq b<m  \} \right )\subseteq \R_{\geq 0}^d.
$$

We point out that $S(f)$ depends on the chosen expansion (1.1). However, it is immediately seen that the
polyhedron $\Delta(f;{\bf u};y)$ (definition below, $f$, $u_1, \ldots ,u_d$ and $y$ being fixed) is independent
of the chosen expansion (1.1): indeed, the semigroup $<{\rm Supp}(f)>\subset \N^{d+1}$ generated by the support of $f$ obviously is.
Of course the $F$-subset $\Delta(f;{\bf u};y)$ is a rational polyhedron because expansion (1.1) is required to be finite.
See Hironaka's theorem [H](4.8) (restated as proposition {\bf I.4} below)
about the characteristic polyhedron.

\bigskip

\noindent {\bf I.2 DEFINITION.}{\it With notations as above:

\noindent (i) the rational polyhedron $\Delta(f;{\bf u};y) \subset \R_{\geq 0}^d$
is the smallest $F$-subset containing all points of $S(f)$,
$f$ expanded as in (1.1), each $C_{{\bf a},b}$ a unit in $R$;

\noindent (ii) the ``characteristic polyhedron'' $\Delta(f;{\bf u}) \subset \R_{\geq 0}^d$ is the $F$-subset
defined by the formula
$$
\Delta(f;{\bf u}):=\bigcap_{(\hat{y},u_1,\dots,u_d)}\Delta(f;{\bf u};\hat{y}), \eqno (1.4)
$$
where the intersection runs over all r.s.p's of $\hat{R}$ of the form $(\hat{y},u_1,\dots,u_d)$;

\noindent (iii) let $L: \ (x_1,x_2,\dots,x_d)\mapsto L(x_1,x_2,\dots,x_d)=\lambda_{1}x_1+\lambda_{2}x_2+\dots+\lambda_{d}x_d$,
$\lambda_{1},\lambda_{2},\dots,\lambda_{d}\in \Q_{\geq 0}$ be a nonzero nonnegative linear form on $\R^d$. We define
$$
l (f,{\bf u},y):={\rm min}\{{L({\bf a}) }|{\bf a}\in \Delta(f;{\bf u};y)\}\geq 0.
$$
We define a {\rm monomial valuation} $v_{L,{\bf u},y,f}$ on $R$ by setting
$$
I_\lambda := (\{y^b{\bf u}^{\bf a}| l(f,{\bf u},y)b+ L({\bf a})\geq \lambda \})\subseteq R,
$$
for $\lambda \geq 0$ and $v_{L,{\bf u},y,f}(g):={\rm min}\{\lambda \in \Q | g \in I_\lambda\}$ for any nonzero $g \in R$.
}

\medskip

The following proposition is an easy but useful exercise left to the reader.

\medskip

\noindent {\bf I.3 PROPOSITION.} {\it Let $L$ be a nonzero nonnegative linear form as above, and let
$$I:=\{i | \lambda_i > 0\}, \ I':=\{i | \lambda_i= 0\}= \{1, \ldots ,d\} \backslash I.
$$
The graded algebra $ \gr_{v_{L,{\bf u},y,f}}(R)$ of $R$ w.r.t. $v_{L,{\bf u},y,f}$ is given by

\noindent (i) if $l(f,{\bf u},y)\not=0$, then
$$
\gr_{v_{L,{\bf u},y,f}}(R)={R \over (y,\{u_i\}_{i\in I})}[Y,\{U_i\}_{i\in I}];
$$

\noindent (ii) if $l(f,{\bf u},y)=0$, then
$$
\gr_{v_{L,{\bf u},y,f}}(R)={R \over (\{u_i\}_{i\in I})}[\{U_i\}_{i\in I}].
$$

\noindent In particular, we have $\gr_{v_{L,{\bf u},y,f}}(R)\simeq k[Y,U_1,U_2,\dots,U_d]$ whenever $L$ is positive.

}

\bigskip

The following is Hironaka's theorem [H](4.8).

\medskip

\noindent {\bf I.4 PROPOSITION (Hironaka).} {\it With notations as above, there exists  $\hat{y} \in \hat{R}$ such that
$(\hat{y},u_1,\dots,u_d)$ is a r.s.p. of $\hat{R}$ and $\Delta(f;{\bf u};\hat{y})= \Delta(f;{\bf u})$. In particular,
$\Delta(f;{\bf u})$ is a rational polyhedron.

}

\bigskip

\centerline {\underbar{II GETTING RATIONAL COORDINATES}}

\medskip

We now introduce the G-ring assumption on $R$ ([EGAIV] 7.3.13 w.r.t. that property {\bf P} in 7.3.8 (iv)'
defined in 6.7.6; see also [M] top of p.256). We briefly review this notion in the special case
of a regular local ring $R$ as above.

\smallskip

Excellent rings are defined in [EGAIV] (7.8.2). This definition consists of three properties:

\smallskip

(i) universal catenariness; (ii) geometrical regularity of formal fibers;
(iii) openness of the regular locus of domains which are finite $R$-algebras.

\smallskip

Rings satisfying only (ii) and (iii) are called quasi-excellent. Rings satisfying only (ii) are called $G$-rings.
Summing up the material which is useful for our purpose, we have:

\medskip

\noindent {\bf II.1 LEMMA.} {\it Let $(R,M,k)$ be a regular local ring. Then:
$$
R \ {\rm is} \ {\rm excellent} \Leftrightarrow R \ {\rm is} \ {\rm quasi-excellent}
\Leftrightarrow R \ {\rm is} \ {\rm a} \ {\rm G-ring},
$$
and this is again equivalent to:
$$
\forall P\in {\rm Spec}R, \ \hat{R}\otimes_R\kappa (P) \ {\rm is} \ {\rm geometrically} \ {\rm regular}. \eqno (2.1)
$$

Assume that a regular local ring $(R,M,k)$ is a G-ring. The following holds:

\smallskip

\noindent (a) let $f \in M$, $f\neq 0$ be such that $R/(f)$ is a domain. Then $\hat{R}/(f)$ is reduced;

\noindent (b) any quotient of $R$, localization $R_P$ of $R$ at a prime $P\in {\rm Spec}R$
or localization of the polynomial ring $R[T]$ at a maximal ideal is again a G-ring.

}

\medskip
\noindent {\it Proof.} For the first statement, we need to check that (2.1) implies (i), (ii) and (iii) in
the definition of excellent rings [EGAIV] (7.8.2) for any regular local ring $R$. In our context,
any regular ring satisfies (i) by [EGAIV] (7.1.10) and (7.1.11)(i). Property (ii) on formal fibers
may be checked only for the maximal ideal $\frak{p}=M$ by {\it ibid.} (7.8.3)(i), which is precisely (2.1).
Also (iii) is a consequence of (ii) for $R$ local, {\it ibid.} (7.8.3)(i). This proves the first statement.

\smallskip

Let $f=\hat{\delta} \hat{f}_1^{a_1}\cdots \hat{f}_s^{a_s}$ be a decomposition of $f$ into pairwise
nonequivalent irreducible factors in the UFD $\hat{R}$, where $\hat{\delta}\in \hat{R}$ is a unit. By (2.1),
the formal fiber ring $\hat{R}\otimes_RQF(R/(f))$ is reduced, so $a_1=\cdots =a_s=1$ and (a) is proved.
Finally (b) is a special case of [EGAIV] (7.8.3)(ii).\hfill{}$\square$

\medskip

\noindent {\bf II.2 LEMMA.} {\it Let $(R,M, k)$ be a regular local ring with r.s.p.
$(y,{\bf u})$ which is a G-ring. Let $f\in M$, $f \not \in (u_1, \ldots ,u_d)$ and assume that
$\Delta (f;{\bf u}) = \emptyset $.

\medskip

Then there exists $z \in R$ such that $(z,u_1,\dots,u_d)$ is a r.s.p. of $R$ and
$$
\Delta (f;{\bf u};z) = \emptyset .
$$
}

\medskip

\noindent {\it Proof.} By Hironaka's theorem (proposition {\bf I.4}), there exists a r.s.p.
$(\hat{y},u_1,\dots,u_d)$ of $\hat{R}$ such that $\Delta (f;{\bf u};\hat{y}) = \emptyset $. By definition {\bf I.2}.(i), this means
that $f=\hat{\delta}\hat{y}^m$ for some unit $\hat{\delta} \in \hat{R}$.

Let $(z) \subset R$ be a prime divisor of $(f)$. By lemma {\bf II.1}(a), $\hat{R} /(z)$ is reduced, so
$(z) \hat{R}=(\hat{y})$. By faithful flatness of completions, $(z)=(\hat{y})\cap R$ and this proves that
$\sqrt{(f)}=(z)$ is prime. We have
$$
{\rm ord}\overline{z}=1,  \ (f)=(z)^m,
$$
i.e. $(z,u_1,\dots,u_d)$ is a r.s.p. of $R$ and $\Delta (f;{\bf u};z) = \emptyset $.\hfill{}$\square$

\bigskip

\noindent {\it Note.} The proof of the main theorem below uses a different algorithm from Hironaka's. This is
illus\-tra\-ted in example {\bf II.5} below.

\bigskip

\noindent {\bf II.3 THEOREM.} {\it Let $(R,M, k)$ be a regular local ring with r.s.p.
$(y,{\bf u})$ which is a G-ring, and $f\in M$, $f \not \in (u_1, \ldots ,u_d)$.

\medskip

Then there exists $z \in R$ such that $(z,u_1,\dots,u_d)$ is a r.s.p. of $R$ and
$$
\Delta (f;{\bf u};z) = \Delta (f;{\bf u}) . \eqno (2.2)
$$
}

\medskip
\noindent {\it Proof.} By lemma {\bf II.2}, it can be assumed that $\Delta (f;{\bf u}) \neq \emptyset$.
By Hironaka's theorem (proposition {\bf I.4}), there exists a r.s.p. of the form $(\hat{y},u_1, \ldots ,u_d)$ of
${\hat R}$ such that $ \Delta(f;{\bf u};\hat{y})=\Delta(f;{\bf u})$.
The rational polyhedron $\Delta (f;{\bf u})$ may then be defined by a formula
$$
\Delta (f;{\bf u})=\{  {\bf x}=(x_1,x_2,\dots,x_d)\in \R^d_{\geq 0} \ \vert \ L_j(x_1,x_2,\dots,x_d) \geq 1,\  1\leq j \leq n \},
$$
for some $n \geq 1$; here each $L_j: \ \R^d \rightarrow \R^d$ is a (nonzero) nonnegative rational linear form verifying
$$
L_j(\Delta (f;{\bf u}))=\lbrack 1,+\infty \lbrack . \eqno (2.3)
$$
We fix $A \geq 1$ such that for all $j$, $1 \leq j \leq n$, we have
$$
L_j(x_1,x_2,\dots,x_d):=\lambda_{j,1}x_1+\lambda_{j,2}x_2+\cdots +\lambda_{j,d}x_d, \
\lambda_{j,1},\lambda_{j,2},\ldots, \lambda_{j,d}\in {1 \over A}\N.
$$

\medskip

\medskip

With notations as in {\bf I.2}.(iii), we let
$$
L_j(\Delta (f;{\bf u};y))=:\lbrack l_j(f,{\bf u},y),\infty \lbrack , \ 1\leq j \leq n .
$$

Note that $L_j(x_1,x_2,\dots,x_d)=l_j(f,{\bf u},y)$ is the equation of a face of  $\Delta(f;{\bf u};y)$. We attach a rational number:
$$
\Lambda (f,{\bf u},y):= \sum_{j=1}^n{(1 -l_j(f,{\bf u},y))}\geq 0. \eqno (2.4)
$$
Following Hironaka [H] (2.6), we consider the initial form ${\rm in}_{v_j}f$ of $f$ with respect to the
valuation $v_j:=v_{L_j,{\bf u},y,f}$ of {\bf I.2}.(iii).

\medskip

Suppose equality (2.2) is not achieved with $z=y$. By (2.3), this means that $l_j(f,{\bf u},y)<1$ for {\it some} $j$,
$1 \leq j \leq n$ which we fix now. In particular we have $\Lambda (f,{\bf u},y)>0$. We consider two cases:

\medskip

\noindent {\it Case 1:} $l_j(f,{\bf u},y)>0$. With notations as in
proposition {\bf I.3}, we let $t_{i'}:={\rm in}_{v_j}u_{i'}\in \gr_{v_j}(R)_0=R/(y,\{u_{i}\}_{i\in I})$ for each $i'\in I'$;
let $t_i:=U_i\in \gr_{v_j}(R)$ for $i\in I$.

\medskip

\noindent {\it Case 2:} $l_j(f,{\bf u},y)=0$. Similarly, we let $t_{i'}:={\rm in}_{v_j}u_{i'}\in
\gr_{v_j}(R)_0=R/(\{u_{i}\}_{i\in I})$ for each $i'\in I'$; let $Y:={\rm in}_{v_j}y\in \gr_{v_j}(R)_0$;
let finally $t_i:=U_i\in \gr_{v_j}(R)$ for $i\in I$.

\medskip

Then
$$
S:= \gr_{v_j}(R)_{(Y,t_1, \ldots ,t_d)}
$$
is a regular local ring with r.s.p. $(Y,t_1, \ldots ,t_d)=:(Y,{\bf t})$ and residue field $S/N=k$,
$N:=(Y,t_1, \ldots ,t_d)$. Note that $S$ is canonically endowed with a monomial
valuation (still denoted by $v_j$) which is induced from the graded structure of $\gr_{v_j}(R)$. The valuation $v_j$
also canonically extends to
the formal completion $\hat{S}$ w.r.t. $N$. We have inclusions $S\subseteq \gr_{v_j}(\hat{R})_{(Y,t_1, \ldots ,t_d)}
\subseteq \hat{S}$, with an isomorphism
$\hat{S}\simeq \gr_{v_j}(\hat{R})_0\lbrack \lbrack Y, \{U_{i}\}_{i\in I} \rbrack \rbrack $ in case 1
(resp. $\hat{S}\simeq \gr_{v_j}(\hat{R})_0\lbrack \lbrack \{U_{i}\}_{i\in I} \rbrack \rbrack $ in case 2). We have
${\rm in}_{v_j}f \in N$, ${\rm in}_{v_j}f \not \in (t_1, \ldots ,t_d)$, so ${\rm in}_{v_j}f \in S$ verifies assumptions
(1.2) and (1.3).

\medskip

Hironaka's construction [H] (vertex preparation, lemma (3.10)) implies the following:

\medskip

\noindent (a) $\hat{Y}:={\rm in}_{v_j}\hat{y}\in \gr_{v_j}(\hat{R})_{l_j(f,{\bf u},y)}$ and $(\hat{Y},t_1, \ldots ,t_d)$ is
a r.s.p. of $\hat{S}$;

\medskip

\noindent (b) $\Delta (\hat{y}; {\bf u};y)\subseteq \Delta (f;{\bf u};y)$.

\medskip

Since $ \Delta(f;{\bf u};\hat{y})=\Delta(f;{\bf u})$, there a finite expansion (1.1) of the form
$$
f=\hat{\gamma}\hat{y}^m +\sum_{0\leq b \leq m-1}\hat{y}^{b}\sum_{{\bf a}} \hat{C}_{{\bf a},b}{\bf u}^{\bf a}, \ {\bf a}\in \N^d,
$$
where $\hat{\gamma}$ and each $\hat{C}_{{\bf a},b}$ is a unit in $\hat{R}$. Let $({\bf a},b)$, $0 \leq b \leq m-1$ appear 
in the above formula. By (2.3), we have $L_j({\bf a})\geq m-b$. Hence
$$
v_j (\hat{C}_{{\bf a},b}\hat{y}^{b}{\bf u}^{\bf a})\geq bl_j(f,{\bf u},y)+m-b> ml_j(f,{\bf u},y)=v_j (\hat{\gamma}\hat{y}^m),
$$
since $(m-b)(1-l_j(f,{\bf u},y))>0$. Computing initial forms, we deduce 

\medskip

\noindent (c) ${\rm in}_{v_j}f=\hat{\delta}\hat{Y}^m$ for some unit 
$\hat{\delta}:={\rm in}_{v_j}\hat{\gamma}\in \gr_{v_j}(\hat{R})_0\subset \hat{S}$, that is: 
$\Delta ({\rm in}_{v_j}f ; {\bf t}; \hat{Y})=\emptyset$.

\medskip

By lemma {\bf II.1}(b), $S$ is a G-ring; by (a), $(\hat{Y},t_1, \ldots ,t_d)$ is
a r.s.p. of $\hat{S}$. Since  $\Delta ({\rm in}_{v_j}f ; {\bf t}; \hat{Y})=\emptyset$ by (c),
lemma {\bf II.2} applies to ${\rm in}_{v_j}f \in N$, ${\rm in}_{v_j}f \not \in (t_1, \ldots ,t_d)$.
So there exists $Z \in N$ such that $(Z,t_1, \ldots ,t_d)$ is a r.s.p. of $S$ and
$$
\Delta ({\rm in}_{v_j}f ; {\bf t}; Z)=\emptyset . \eqno (2.5)
$$
We pick for $Z$ a finite expression in the form (1.1):
$$
Z= \sum_{{\bf a},b} C_{{\bf a},b}Y^b {\bf t}^{\bf a}, \ b\in \N,\ {\bf a}\in \N^d,
$$
where each $C_{{\bf a},b} \in S$ is a unit. Since ${\rm in}_{v_j}f \in \gr_{v_j}(R)$ is homogeneous for $v_j$,
it can assumed w.l.o.g. that $Z\in \gr_{v_j}(R)$ and $Z$ is homogeneous (of degree $l_j(f,{\bf u},y)$);
in particular $C_{{\bf a},b} \in \gr_{v_j}(R)_0$ for each value of $(b, {\bf a})$.

Again by (c), $(Z)=(\hat{Y})$ in $\hat{S}$, which in turn implies that
$$
\Delta (Z; {\bf t};Y)=\Delta (\hat{Y}; {\bf t};Y). \eqno (2.6)
$$
Since $\Delta (\hat{Y}; {\bf t};Y) \subseteq \Delta (\hat{y}; {\bf u};y)$, any set of preimages $\gamma_{{\bf a},b} \in R$
of the elements  $C_{{\bf a},b}\in \gr_{v_j}(R)_0$ defines an element
$$
z:=\sum_{{\bf a},b} \gamma_{{\bf a},b}y^b {\bf u}^{\bf a} \in R.
$$
Since for each value of $(b, {\bf a})$, $C_{{\bf a},b}$ is a unit, $\gamma_{{\bf a},b}\in R$ is a unit.
This implies that $\Delta (z; {\bf u};y)=\Delta (Z; {\bf t};Y)$, wherefrom we conclude using (2.6) and (b) that
$$
\Delta (z; {\bf u};y)\subseteq \Delta (f;{\bf u};y). \eqno (2.7)
$$
Now (2.7) implies that $\Delta (f;{\bf u};z)\subseteq \Delta (f;{\bf u};y)$, and in particular we get
$$
l_k(f,{\bf u},z) \geq l_k(f,{\bf u},y), \ 1 \leq k \leq n. \eqno (2.8)
$$
On the other hand (2.5) is equivalent to
$$
l_j(f,{\bf u},z)>l_j(f,{\bf u},y)
$$
and we conclude from (2.4) that
$$
0 \leq \Lambda (f,{\bf u},z) < \Lambda (f,{\bf u},y). \eqno (2.9)
$$

\medskip

Since $\Lambda (f,{\bf u},y') \in {1 \over m!A}\N$ for any r.s.p. $(y',{\bf u})$ of $R$, (2.9) cannot repeat infinitely many times,
that is, $\Lambda (f,{\bf u},z)=0$ for some $z$ and (2.2) eventually holds. \hfill{}$\square$

\bigskip

The following monic version of theorem {\bf II.3} is useful for applications and used in [CP3] theorem 2.4. The conclusion
of corollary {\bf II.4} is of course not true in general for polynomials which are not monic unless $S$ is Henselian
(this is easily seen by taking $m=1$).

\medskip

\noindent {\bf II.4 COROLLARY.} {\it Assume that $R=S[y]_{(N,y)}$, $(S,N)$ a regular local ring with r.s.p.
$(u_1,\dots,u_d)$ which is a G-ring, and
$$
f =y^m +f_1y^{m-1} +\cdots +f_m \in S[t], \ f_i \in N, \ 1 \leq i \leq m.
$$
Then there exists $z :=y - \phi$, $\phi \in N$ such that $\Delta (f;{\bf u};z) = \Delta (f;{\bf u})$ .
}

\medskip

\noindent {\it Proof.} We consider only those expansions (1.1) with $C_{{\bf a},b}\in \N^d\times\{0, 1, \ldots ,m-1\}$
and $C_{{\bf 0},m}=1$, i.e. monic expansions of degree $m$ in $y$.  These expansions are preserved by coordinate changes
of the form $z:=z- \phi$, $\phi \in S$, and the proof is identical. \hfill{}$\square$

\medskip

We illustrate our method with the following simple example.

\medskip

\noindent {\bf II.5 EXAMPLE.} {\it Take $d=2$ and let
$$
f:= (y - \phi)^m +u_1^mu_2, \ \phi \in (u_1,u_2).
$$
Let $\overline{\phi}_i$ be the image of $\phi$ in the discrete valuation ring $R/(y,u_i)$, $i=1,2$, and
assume furthermore that $0<m_i:={\rm ord}\overline{\phi}_i<+\infty$, $i=1$ and $i=2$}.

\medskip

The polyhedron $\Delta (f;u_1,u_2;y)\subset \R^2_{\geq 0}$ has three vertices:
$$
{\bf v}_0:=(1,{1 \over m}), \ {\bf v}_1:=(m_2,0), \ \ {\bf v}_2:=(0,m_1),
$$
where ${\bf v}_0$ is not solvable, and ${\bf v}_1, {\bf v}_2$ are solvable. Taking $z:=y-\phi$ computes
the characteristic polyhedron
$$
\Delta (f;u_1,u_2)=\Delta (f;u_1,u_2;z)={\bf v}_0 +\R^2_{\geq 0}.
$$
Our algorithm first expresses
$$
\Delta (f;u_1,u_2)=\{(x_1,x_2)\in \R^2_{\geq 0}: x_1\geq 1, \ mx_2\geq 1\}.
$$
We thus take $n=2$, $L_1(x_1,x_2):=x_1$, $L_2(x_1,x_2):=mx_2$ in (2.3), so
$$
l_j(f;u_1,u_2;y)=0, \ {\rm in}_{v_j}f=(\overline{y} -\tilde{\phi}_j)^m \in R/(u_j), \ j=1,2.
$$
By assumption on $\phi$, the image of $\tilde{\phi}_j$ in $R/(y,u_j)$ is $\overline{\phi}_i$ which is nonzero.

\smallskip

Taking for instance $j=1$ as a first step of the algorithm, $\gamma_1\in R$ any preimage of the unit
$\overline{u}_2^{-m_1}\overline{\phi}_1 \in R/(y,u_1)$, we let $z_1:=y - \gamma_1u_2^{m_1}$. The
polyhedron $\Delta (f;u_1,u_2;z_1)\subset \R^2_{\geq 0}$ now has two vertices: ${\bf v}_0$ and ${\bf v}_1$.
A second step of the algorithm similarly produces $z_2:=z_1 - \gamma_2u_1^{m_2}$, $\gamma_2 \in R$ a unit
and $\Delta (f;u_1,u_2;z_2)=\Delta (f;u_1,u_2)$.

On the other hand, Hironaka's algorithm only takes {\it vertices} into account and not whole faces as does ours.
It thus introduces a succession of coordinate changes
$$
z(0):=y, \ z(i+1):=z (i) - \gamma (i)u_{j(i)}^{m(i)}, \ \gamma (i) \in R \ {\rm a} \ {\rm unit},
$$
where $(i'>i, j(i')=j(i)) \Longrightarrow m(i')>m(i)$. At each step $j(i)\in \{1,2\}$ is chosen in such a way that
$m(i)$ is minimal among the possible choice for $j(i)$. This ensures that $m(i)$ goes to infinity with $i$ and one gets
convergence in the formal completion
$$
z(\infty) := y -\hat{\phi}, \ \hat{\phi}:=\sum_{i=1}^{\infty}\gamma (i)u_{j(i)}^{m(i)}\in \hat{R}
$$
with no control on $\hat{\phi}$. \hfill{}$\square$

\medskip

The following example in mixed characteristic is a (very) special case of the constructions performed in
the proof of [CP3] corollary 1.2. We compute the $\Delta$-function (0.2) along a permissible curve $Z={\cal C}$
for an arithmetical threefold ${\cal X}$ and sketch the application to Resolution.

\bigskip

\noindent {\bf II.6 EXAMPLE.} {\it Take $A:=\Z [u_2,u_3,y]$, $m \in \Z$, $m \not \equiv 0 \ {\rm mod}p$, and let
$$
f:=y^p -u_2u_3^p -(m +u_3^p)p^p,  \ {\cal X}:={\rm Spec}(A/(f)).
$$
}

Computing ${\rm Disc}_y(f)=p^p\left (u_2u_3^p +(m +u_3^p)p^p\right )^{p-1}$, it is seen that the locus ${\rm Sing}_p{\cal X}$
of multiplicity $p$ of ${\cal X}$ is the curve ${\cal C}:=(y,p,u_3)$.
Let $O:=(p,u_2,u_3,y)\in {\cal C}$, $R:=A_{(p,u_2,u_3,y)}$. In corollary {\bf II.4}, we may take here
$z:=y-(m+u_3)p$. There is an expansion
$$
f=z^p +\sum_{i=1}^{p}f_iz^{p-i} \eqno (2.10)
$$
where
$$
f_i:=\pmatrix {p \cr i \cr}(m+u_3)^ip^i, \ 1 \leq i \leq p-1;
\ f_p:=u_2u_3^p +\left ( m^p-m  + \sum_{i=1}^{p-1}\pmatrix {p \cr i \cr}u_3^im^{p-i} \right )p^p.
$$
Define ${\bf v}_0:=(0,1/p,1)$, ${\bf v}_1:=(1+1/p,0,0)$, ${\bf v}'_1:=(1+1/(p-1),0,0)$, ${\bf v}'_2:=(1+1/p,0,1/p)$.
Using elementary arithmetics, it is checked that the vertices of $\Delta (f;p,u_2,u_3;z)\subset \R^3_{\geq 0}$ are:
$$
\left\{
\matrix{
  {\bf v}_0, {\bf v}_1 \hfill{}& {\rm if} &  m^p-m \not \equiv 0 \ {\rm mod}p^2\cr
   &  &  \cr
  {\bf v}_0, {\bf v}'_1, {\bf v}'_2  & {\rm if} & m^p-m  \equiv 0 \ {\rm mod}p^2 \cr
}
\right.
.
$$
The only vertex in $\N^3$ in the above list is ${\bf v}'_1$ when $p=2$. In this special case, the corresponding
initial form is
$$
{\rm in}_{{\bf v}'_1}(f)=Z^2+P^2Z+\lambda P^4 \in \F_2[P,U_2,U_3,Z], \ \lambda \in \F_2, \ P:={\rm in}p.
$$
Since this polynomial is not a square, this proves that $\Delta (f;p,u_2,u_3;z)=\Delta (f;p,u_2,u_3)$.

\smallskip
Equation (2.10) is defined at every point $P \in {\rm Spec}A$ such that $z \in A_P$. By theorem {\bf II.3},
this always holds on a nonempty Zariski neighborhood $U$ of $O$ ($U={\rm Spec}A$ in this example).
Take a closed point $x \in {\rm Sing}_p{\cal X}$, $x \neq O$. Then $A_x$ has a r.s.p. $(p,w_x:=u_2-\gamma^p, u_3,z)$,
$\gamma \in A_x$ a unit. Taking $z_x:=z -\gamma u_3$, one lets ${\bf v}:=(1/(p-1),0,1)$ and checks that
$$
\Delta (f;p,w_x,u_3; z_x)=\Delta (f;p,w_x,u_3)=
{\rm Conv}\left ({\bf v}+\R^3_{\geq 0}, \Delta (f;p,u_2,u_3)\right )\subset \R^3_{\geq 0}.
$$
Analyzing the effects of blowing up ${\cal C}$ on polyhedra, it can be proved that the blowing up
${\cal X}'$ of ${\cal X}$ along ${\cal C}$ satisfies
$$
\left\{
\matrix{
  {\cal X}' \ {\rm regular} \hfill{}& {\rm if} &  m^p-m \not \equiv 0 \ {\rm mod}p^2 \hfill{}\cr
   {\rm Sing}_p {\cal X}'=\emptyset & {\rm if} & p\neq 2 \ {\rm and} \  m^p-m  \equiv 0 \ {\rm mod}p^2\cr
  {\rm Sing}_2 {\cal X}'= {\cal C}' & {\rm if} & p= 2 \ {\rm and} \  m  \equiv 1 \ {\rm mod}4 \hfill{}\cr
}
\right.
,
$$
where ${\cal C}'$ maps isomorphically to ${\cal C}$ (blowing up ${\cal C}'$ resolves the singularities of
${\cal X}'$ in this last case).

\bigskip

\centerline {\underbar{BIBLIOGRAPHY}}

\bigskip

\noindent [C] COSSART Vincent, Sur le poly\`edre caract\'eristique d'une singularit\'e, {\it Bull. Soc. Math.} {\bf 103}-1 (1975), 13-19.

\medskip

\noindent [CGO] {COSSART Vincent, GIRAUD Jean, ORBANZ Ulrich}, Resolution of surface singularities. With an appendix by H. Hironaka.
{\it Lect. Notes in Math.} {\bf 1101}, Springer-Verlag (1984).

\medskip

\noindent [CP1,2] COSSART Vincent, PILTANT Olivier,
Resolution of singularities of threefolds in positive characteristic I, {\it J. Algebra}
{\bf 320}, no. 7 (2008), 1051-1082; and II {\it ibid.} {\bf 321}, no. 1 (2009) 1836-1976.

\medskip

\noindent [CP3] COSSART Vincent, PILTANT Olivier, Resolution of singularities of arithmetical threefolds I, {\it preprint}
http://hal.archives-ouvertes.fr [hal-00873967].

\medskip

\noindent [EGAIV] GROTHENDIECK Alexander, DIEUDONN\'E Jean,
\'El\'ements de g\'eom\'etrie alg\'ebrique IV, {\it Publ. Math. I.H.E.S.} {\bf 24} (1965).

\medskip

\noindent [H] HIRONAKA Heisuke,  Characteristic polyhedra of  singularities,
{\it J. Math. Kyoto Univ.} {\bf 7}-3 (1967), 251-293.

\medskip

\noindent  [M] MATSUMURA Hideyuki, Commutative ring theory, 3rd edition, {\it Cambridge studies in advanced mathematics} {\bf 8} (1986), Cambridge Univ. Press.

\bye